# BROWNIAN MOVING AVERAGES HAVE CONDITIONAL FULL SUPPORT

By Alexander Cherny

*Moscow State University*

We prove that any Brownian moving average

$$X_t = \int_{-\infty}^{t} (f(s-t) - f(s)) \, dB_s, \qquad t \geq 0,$$

satisfies the *conditional full support* condition introduced by Guasoni, Rásonyi and Schachermayer [*Ann. Appl. Probab.* **18** (2008) 491–520].

**1. Introduction.**

1.1. *Overview.* It is well known (see Soner, Shreve and Cvitanić [8], Levental and Skorokhod [6], Cherny [2]) that in the Black–Scholes–Merton model with proportional transaction costs the superreplication price of a European call option is equal to its trivial upper bound. The same is true for any European type contingent claim in this model (see Cvitanić, Pham and Touzi [3]). In the recent paper [4], Guasoni, Rásonyi and Schachermayer proved that the same result holds for a much wider class of models satisfying only a minor geometric condition termed *conditional full support* and denoted CFS for brevity (see the paper by Kabanov and Stricker [5] for further research in this direction).

The CFS condition is as follows. We consider a filtered probability space $(\Omega, \mathcal{F}, (\mathcal{F}_t)_{t \in [0,T]}, \mathsf{P})$ and a continuous $(\mathcal{F}_t)$-adapted process $(X_t)_{t \in [0,T]}$ meaning the discounted price (or the logarithm of the discounted price) of an asset. The CFS condition requires that, for any $t \in [0, T]$,

$$\operatorname{supp} \operatorname{Law}(X_u; t \leq u \leq T \mid \mathcal{F}_t) = C_{X_t}[t, T] \qquad \text{a.s.,}$$

where $C_x[t, T]$ denotes the space of continuous real-valued functions on $[t, T]$ with $f(t) = x$ and "supp" denotes the support (the conditional distribution

Received July 2007; revised November 2007.

here is viewed as a measure on the space $C[t,T]$ of continuous functions on $[t,T]$).[1]

1.2. *Goal of the paper.* As motivated by the above discussion, the CFS condition is interesting and important. The paper [4] provides several examples of processes satisfying this condition. One of them is the fractional Brownian motion (FBM). It is well known (see Mandelbrot and Van Ness [7]) that FBM is a Brownian moving average, that is, it can be represented as

$$(1.1) \qquad X_t = \int_{-\infty}^{t} (f(s-t) - f(s))\, dB_s, \qquad t \in [0,T],$$

with a certain function $f : \mathbb{R} \to \mathbb{R}$ such that $f = 0$ on $\mathbb{R}_+$ and $\int_{-\infty}^{t}(f(s-t) - f(s))^2\, ds < \infty$ for any $t \geq 0$. Let us remark that the class of moving averages includes processes that are, in a sense, more convenient for financial modeling than FBM; for example, FBM is not a semimartingale (except for two particular cases), while a moving average is a semimartingale provided that $f$ is absolutely continuous and its derivative is square integrable on $(-\infty, 0]$ (see Cheridito [1]).

The main result of the paper is

THEOREM 1.1. *Let $f : \mathbb{R} \to \mathbb{R}$ be a function such that $f = 0$ on $\mathbb{R}_+$, $\int_{-\infty}^{t}(f(s-t) - f(s))^2\, ds < \infty$ for any $t \geq 0$, and $f$ is not zero on a set of positive Lebesgue measure. Then the process $X$ defined by (1.1) satisfies the CFS condition with respect to its natural filtration.*

We also consider the CFS condition for general Gaussian processes. In discrete time it is easy to see that the CFS condition (appropriately redefined for the discrete-time case) is satisfied provided that $X$ is a Gaussian process such that $\mathrm{Var}(X_t - X_s \mid X_u; u \leq s) > 0$ for any $s < t$ (by Var we denote the variance). This might seem a bit surprising, but in continuous time the corresponding result does not hold; see Example 3.1.

**2. Proof of Theorem 1.1.** Let $T > 0$ and let $f \in L^2[-T, 0]$. For $g \in L^2[0,T]$, we denote by $f * g$ the convolution of $f$ and $g$ restricted to $[0,T]$, that is, the function

$$(f * g)(t) = \int_0^t f(s-t) g(s)\, ds, \qquad t \in [0,T].$$

---

[1] We deal with real-valued price processes, while Guasoni, Rásonyi and Schachermayer deal with strictly positive processes. The relationship between the two definitions is trivial: a process $X$ satisfies our version of CFS if and only if $e^X$ satisfies the CFS condition from [4].

LEMMA 2.1. *Let $h \in L^2[-T, 0]$ satisfy the condition $\int_{-\varepsilon}^{0} |h(t)| \, dt > 0$ for any $\varepsilon > 0$. Then the space $\{h * g : g \in L^2[0, T]\}$ is dense in $C_0[0, T]$.*

PROOF. If $g$ is absolutely continuous with a square-integrable derivative and $g(0) = 0$, then $(h * g)' = h * g'$. Thus, if a function $h * g$ approximates a function $\varphi \in L^2[0, T]$ in the $L^2$-sense, then the function $h * G$, where $G(t) = \int_0^t g(s) \, ds$, approximates the function $\Phi(t) = \int_0^t \varphi(s) \, ds$ in the $C_0[0, T]$-sense. So, it is sufficient to prove that the space $\{h * g : g \in L^2[0, T]\}$ is dense in $L^2[0, T]$.

Suppose that this is not true. Then there exists a function $\varphi \in L^2[0, T]$ not identically equal to zero and such that

$$\int_0^T (h * g)(t) \varphi(t) \, dt = 0 \qquad \forall g \in L^2[0, T].$$

This means that

$$0 = \int_0^T \int_0^t h(s - t) g(s) \varphi(t) \, ds \, dt$$

$$= \int_0^T \int_s^T h(s - t) g(s) \varphi(t) \, dt \, ds \qquad \forall g \in L^2[0, T],$$

which, in turn, is equivalent to the property

$$\int_s^T h(s - t) \varphi(t) \, dt = 0 \qquad \forall s \in [0, T].$$

But this is impossible due to the Titchmarsh convolution theorem (see [9], Chapter VI). The obtained contradiction yields the desired result. □

PROOF OF THEOREM 1.1. Let $a \in (-\infty, 0]$ be a number such that $f = 0$ a.e. with respect to the Lebesgue measure on $[a, 0]$ and $\int_{a-\varepsilon}^{a} |f(x)| \, dx > 0$ for any $\varepsilon > 0$. We can assume that $a = 0$. The case $a < 0$ is reduced to this one by considering the new Brownian motion $\widetilde{B}_t = B_{t-a} - B_{-a}$ and the new function $\widetilde{f}(x) = f(x - a)$.

We have to prove that, for any $t \in [0, T]$,

$$\operatorname{supp} \operatorname{Law}(X_u - X_t; t \leq u \leq T \mid \mathcal{F}_t) = C_0[t, T] \qquad \text{a.s.},$$

where $\mathcal{F}_t = \sigma(X_s; s \leq t)$. Obviously, it is sufficient to prove the above property with $\mathcal{F}_t$ replaced by the larger filtration $\mathcal{G}_t = \sigma(B_s : -\infty < s \leq t)$. With this substitution, it is obviously sufficient to check the property only for $t = 0$. We then have

$$\operatorname{Law}(X_u; 0 \leq u \leq T \mid \mathcal{G}_0)(\omega)$$

$$= \operatorname{Law}\left(\int_0^u f(v - u) \, dB_v\right.$$

$$+ \int_{-\infty}^{0} (f(v-u) - f(v)) \, dB_v; 0 \leq u \leq T \, \Big| \, \mathcal{G}_0 \Big)(\omega)$$

$$= \mathrm{Law}\bigg( \int_{0}^{u} f(v-u) \, dB_v + \varphi(u, \omega); 0 \leq u \leq T \bigg),$$

where $\varphi(\cdot, \omega)$ is the path of the process $Y = \int_{-\infty}^{0} (f(v - \cdot) - f(v)) \, dB_v$ corresponding to the elementary outcome $\omega$.

The above equality means that the conditional law of $(X_u)_{u \in [0,T]}$ given $\mathcal{G}_0$ is nothing but the unconditional law of $(\int_{0}^{u} f(v-u) \, dB_v)_{u \in [0,T]}$ shifted by the function $\varphi(u, \omega)$. As the two laws differ by such a shift, it is sufficient to prove that

$$(2.1) \qquad \mathrm{supp} \, \mathrm{Law}\bigg( \int_{0}^{u} f(v-u) \, dB_v; 0 \leq u \leq T \bigg) = C_0[0, T].$$

It follows from the Girsanov theorem that, for any $g \in L^2[0, T]$,

$$\mathrm{Law}\bigg( \int_{0}^{u} f(v-u) \, dB_v; u \leq T \bigg)$$

$$\sim \mathrm{Law}\bigg( \int_{0}^{u} f(v-u) \, dB_v + \int_{0}^{u} f(v-u) g(v) \, dv; u \leq T \bigg).$$

Hence, if a function $\psi$ belongs to the left-hand side of (2.1), then the same is true for $\psi + \int_{0}^{\cdot} f(v - \cdot) g(v) \, dv$. Using now the nonemptiness of the support and recalling Lemma 2.1, we obtain (2.1), which completes the proof. □

**3. Example.** Let $(X_n)_{n=0,\ldots,N}$ be a Gaussian random sequence such that

$$(3.1) \qquad \mathrm{Var}(X_n - X_{n-1} \mid X_i; i \leq n-1) > 0 \qquad \forall n = 1, \ldots, N.$$

Using induction in $m$, it is then easy to see that $X$ satisfies the discrete-time version of the CFS condition:

$$(3.2) \qquad \mathrm{supp} \, \mathrm{Law}(X_i : i = n+1, \ldots, m \mid X_i : i = 0, \ldots, n) = \mathbb{R}^{m-n}$$

$$\forall 0 \leq n < m \leq N.$$

Let us remark that (3.2) obviously implies (3.1), so that the latter property serves as a criterion for the CFS for discrete-time Gaussian processes.

Surprisingly enough, in continuous time such a simple criterion does not hold, as shown by the next example.

EXAMPLE 3.1. Let $B$ be a Brownian motion. For $n \in \mathbb{Z}_+$, denote $a_n = 1 - 2^{-n}$ and let

$$X_t^n = b_n \int_0^t I(a_n \leq s \leq a_{n+1}) \, dB_s$$

$$+ b_n 2^{2n+3} \int_{a_n}^{1} (B_{s \wedge a_{n+1}} - B_{a_n}) \, ds \int_0^t I(s \geq a_{n+1}) \, ds, \qquad t \in [0, 1].$$

The constants $b_n$ are strictly positive and decrease to zero fast enough to ensure that

$$\sum_{n=0}^{\infty} \sup_{t \in [0,1]} |X_t^n| < \infty \qquad \text{a.s.}$$

Then the process

$$X_t = \sum_{n=0}^{\infty} X_t^n, \qquad t \in [0,1]$$

is continuous and Gaussian. For any $0 \leq s < t \leq 1$, the difference $X_t - X_s$ can be represented as $\xi_1 + \xi_2$, where $\xi_1$ is $\sigma(X_u; u \leq s)$-measurable and $\xi_2$ is nondegenerate and depends on the increments of $B$ after time $s$. Hence,

$$\operatorname{Var}(X_t - X_s \mid X_u; u \leq s) > 0 \qquad \forall 0 \leq s < t \leq 1.$$

On the other hand,

$$\int_0^1 X_t \, dt = \sum_{n=0}^{\infty} \int_0^1 X_t^n \, dt$$
$$= \sum_{n=0}^{\infty} b_n \int_{a_n}^1 (B_{s \wedge a_{n+1}} - B_{a_n}) \, ds \left[ 1 + 2^{2n+3} \int_{a_{n+1}}^1 (s - a_{n+1}) \, ds \right] = 0,$$

so that the CFS condition is violated for $X$ already for $t = 0$.

**Acknowledgments.** This paper would not have been written without the advice of Stanislav Molchanov, who suggested using the Titchmarsh convolution theorem, which is at the heart of the proof of the main theorem of the paper. I express my thanks to Martin Schweizer for having attracted my attention to the paper by Guasoni, Rásonyi and Schachermayer. I am thankful to Youri Kabanov and Miklos Rásonyi for the careful reading of the manuscript and important remarks. I express my thanks to two anonymous referees and an Associate Editor for a number of suggestions that improved the quality of the paper.

DEPARTMENT OF PROBABILITY THEORY
FACULTY OF MECHANICS AND MATHEMATICS
MOSCOW STATE UNIVERSITY
119992 MOSCOW
RUSSIA
E-MAIL: alexander.cherny@gmail.com